\documentclass[12pt]{article}
\usepackage{amssymb,amsmath,amsfonts}
\newtheorem{theorem}{Theorem}
\newtheorem{lemma}{Lemma}
\newtheorem{proposition}{Proposition}
\newtheorem{corollary}{Corollary}
\newtheorem{definition}{Definition}

\begin{document}

$\;$

\begin{center}

\LARGE{ Morse functions on the moduli space\\
\noindent of $G_2$ structures  }

\vspace{1pc}

\large{Sung Ho Wang \\
 Department of Mathematics \\
 Postech \\
 Pohang, Korea 790 - 784 \\
 \emph{email}: \textbf{wang@postech.ac.kr}}
\end{center}

\vspace{2pc}

\begin{center}
\textbf{Abstract}
\end{center}
Let $\, \mathfrak{M}$ be the moduli space of torsion free
$\, G_2$ structures on a compact 7-manifold $\, M$, and let
$\, \mathfrak{M}_1 \subset \mathfrak{M}$ be the $\, G_2$ structures
with volume($M$) $=1$. The cohomology map $\, \pi^3: \mathfrak{M} \to H^3(M, R)$
is known to be a local diffeomorphism. It is proved that every nonzero
element of $\, H^4(M, R) = H^3(M, R)^*$ is a Morse function on
$\, \mathfrak{M}_1 $ when composed with $\, \pi^3$. When dim $H^3(M, R) = 2$, the result in particular implies
$\, \pi^3$ is one to one on each connected component of $\, \mathfrak{M}$.
Considering the first Pontryagin class $\, p_1(M) \in H^4(M, R)$, we formulate
a compactness conjecture on the set of $\, G_2$ structures of volume($M$) $=1$
with bounded $L^2$ norm of curvature, which would imply that every connected
component of $\, \mathfrak{M}$ is contractible. We also observe the locus
$\, \pi^3(\mathfrak{M}_1) \subset H^3(M, R)$ is a hyperbolic affine sphere
if the volume of the torus $\, H^3(M, R) / H^3(M, Z)$ is constant on
$\, \mathfrak{M}_1$.

\vspace{1pc}


\noindent \textbf{MS classification}: 53C38

\thispagestyle{empty}

\newpage
\setcounter{page}{1}

$\;$

\begin{center}
\LARGE{ Morse functions on the moduli space\\
\noindent of $G_2$ structures }

\vspace{1pc}

\large{Sung Ho Wang}
\end{center}

\vspace{1pc}

\section{Introduction}
The moduli space of complex structures on a compact Riemann
surface of genus $1$ or $\geq 2$ can be identified with the
deformation space of Riemannian metrics of constant curvature
$0$(up to scale) or $-1$ respectively, while the latter definition
naturally gives rise to the Weil-Petersson metric.

Let $M$ be a compact, oriented, and spin manifold of dimension 7.
Then $M$ admits a differential 3-form $\phi$ of generic type
called a \emph{definite(positive)} 3-form(\textbf{Section 2})[Br],
and such $\phi$ determines a unique Riemannian metric
$g_{\phi}$ and an orientation on $M$. $\phi$ is called a
\emph{$G_2$ form} if $d\phi=0$, $d *_{\phi} \phi =0$, and the
orientation determined by $\phi$ agrees with the given orientation
of $M$, where $*_{\phi}$ is the Hodge star operator with respect
to $g_{\phi}$. The stabilizer of a definite 3-form in the
Euclidean space $R^7$  is isomorphic to compact simple Lie group
$G_2$, and thus the existence of a $G_2$ form is equivalent to a
torsion free $G_2$ structure on $M$. Throughout this paper,
a $G_2$ structure would mean a torsion free $G_2$ structure and
a $G_2$ manifold would mean a manifold with a $G_2$ structure.
Note the holonomy of the underlying metric of
a $G_2$ manifold is necessarily a subgroup of $G_2$.

It is known that the moduli space of $G_2$ structures, denoted by
$\mathfrak{M}$, is a smooth manifold of dimension $b^3 =$ dim
$H^3(M, R)$. When $M$ has full holonomy $G_2$, or equivalently
$b^1 =$ dim $H^1(M, R)$= 0, a connected component of
$\mathfrak{M}$ coincides with the (Ricci flat) Einstein
deformation space of the underlying $G_2$ metric (the property for
a Ricci flat metric to support a parallel spinor is preserved
under Einstein deformation). In this perspective, one of the
motivation for the present work is to examine the analogy or the
difference between classical Teichm\"uller theory and the
deformation of $G_2$ structures [Tr].

Another motivation comes from the question : \emph{ Can one find
the best $G_2$ form (structure) on a given $G_2$ manifold ?} A
natural condition would be to require a $G_2$ form $\phi$ to
satisfy $[*_{\phi} \phi] = - p_1(M)$, where $p_1(M) \in H^4(M, R)$
is (nonzero) Pontryagin class. Such $\phi$ is best in the sense
that it locally minimizes the $L^2$ norm of the curvature of the
associated metric $g_{\phi}$ among the set of $G_2$ forms with
fixed volume $Vol_{\phi}(M)$.

The main result of this paper is that  any nonzero element $\beta
\in H^4(M, R) = H^3(M, R)^*$ composed with the cohomology map
$\mathfrak{M} \to H^3(M, R)$ is a Morse function on
$\mathfrak{M}_1$, where $\mathfrak{M}_1 \subset \mathfrak{M}$ is
the moduli space of $G_2$ structures of volume 1. In particular,
when $b^1=0$ and $\beta = - p_1(M) \neq 0$, every critical point
is a positive local minimum. As a corollary, $p_1(M) \neq 0$
implies each connected component of $\mathfrak{M}_1$ is
noncompact. We also prove Torelli theorem in case $b^1=0$ and
$b^3=2$ : the cohomology map $\mathfrak{M} \to H^3(M, R)$ is one
to one on each connected component.

For further application of this result toward the question raised
above, an analogue of Mumford compactness theorem for Riemann
surfaces would play a role. Based on the geometry of the Torelli
map(to be defined below) and the main results above, we formulate
a conjecture on the compactness of the set of $G_2$ structures of
volume 1 with uniformly bounded height with respect to $-p_1(M)$.
Specifically, let $\{\phi_n\}$ be a sequence of $G_2$ forms of
volume 1 in a connected component of the moduli space of
$G_2$ forms on a compact manifold $M$ with $b^{1}=0$. Suppose
$\{ - p_1(M)([\phi_n]) \}$ is bounded from above. Does there exists a
subsequence $\{ \phi_{n_k} \}$, a sequence of diffeomorphisms
$f_k$ isotopic to the identity,  and a $G_2$ form $\phi$ of volume
1 such that $f_k^*  \phi_{n_k} \, \to \, \phi$ in $C^1$ ?  If this
compactness theorem is true, a simple Morse theory argument then
implies that each connected component $\mathfrak{M}$ is
contractible, and that a $G_2$ manifold $M$ supports a $G_2$ form
$\phi$ unique up to diffeomorphism such that  $[*_{\phi} \phi] = -
p_1(M)$.

Torelli theorem for compact Riemann surfaces states that the
period map of a Riemann surface determines its complex structure
[Gr]. An essential characterization of the image of the period
map, inversion problem, is also known [AM]. In the last section,
we consider the Torelli map $\pi : \, \mathfrak{M} \to H^3(M,R)
\oplus H^4(M, R)$ defined by
\begin{equation}
\pi(\langle \phi \rangle) = (\pi^3(\langle \phi \rangle), \,
\pi^4(\langle \phi \rangle) )
= ([\phi], [*_{\phi}\phi])
\in H^3(M, R) \oplus H^4(M, R) \notag
\end{equation}
where $\langle \phi \rangle$ is the equivalence class of $G_2$
forms represented by $\phi$. It is proved that $\pi$ is a
Lagrangian immersion, and that it is an isometric immersion(up to
sign) along $\mathfrak{M}_1$ when $b^{1}=0$, where $H^3(M, R)
\oplus H^4(M, R)$ is equipped with  the canonical symplectic form
and  the metric of signature ($b^3, \, b^4$). If the ratio of
Jacobian of the projections $\pi^3$ and $\pi^4$ is constant along
$\mathfrak{M}_1$, each hypersurface $\pi^3(\mathfrak{M}_1) \subset
H^3(M, R)$ and $\pi^4(\mathfrak{M}_1) \subset H^4(M, R)$ describes
a hyperbolic affine sphere centered at the origin. This is also
equivalent to the volume of the torus $ \, H^{3}(M, R) / H^{3}(M, Z)$
being constant along $\mathfrak{M}_1$ [Hi].

\section{Moduli space of $G_2$ structures}
Consider an exterior 3-form on $R^7$
\begin{equation}
\varphi = ( dx^{12} + dx^{34} + dx^{56} ) \wedge dx^7 + dx^{135}
  - dx^{146} - dx^{362} - dx^{524},
\end{equation}
where $x^i$'s are coordinate functions and $dx^{12} = dx^1 \wedge
dx^2$, etc. The stabilizer of
$\varphi$ is a compact simple simply connected Lie group
$G_2  \subset GL(7, R)$, and it also preserves the 7-form $dx^{12..7}$
and a metric uniquely determined by
\begin{equation}
\langle \,  u, v \,\rangle \, dx^{12..7} \; = \; \frac{1}{6} \,( u \,
\lrcorner \, \varphi ) \,
\wedge \, ( v \, \lrcorner \, \varphi ) \, \wedge \, \varphi
\end{equation}
for $u, \, v \, \in R^7$ [Br].
The orbit of $\varphi$ is one of the two open orbits in $\bigwedge^3
(R^7)^*$ under $GL(7, R)$ action.

Let $M$ be a compact oriented 7-manifold. Let
$PM \subset \bigwedge^3 T^*M$ be the open subset whose fiber at each
point $p \in M$ consists of $\phi_p \in \bigwedge^3 T_p^*M$ that can be
identified with (1) under an oriented isomorphism between $T_p M$ and
$R^7$. An element $\phi \in C^{\infty}(PM)$ is called a
\emph{positive(definite)} 3-form. $M$ admits a positive 3-form if and only if
it is spin, and a positive 3-form $\phi$ in turn determines a metric
$g_{\phi}$ and a volume form $\mbox{dvol}_{\phi}$ by (2) [Br]. We denote
by $*_{\phi}$ the Hodge star operator on differential forms defined by
$g_{\phi}$. Since $PM \subset \bigwedge^3 T^*M$ is an open subset,
the tangent space $T_{\phi} C^{\infty}(PM)$ can be
identified with  the space of differential 3-forms
\begin{equation}
T_{\phi} C^{\infty}(PM)  \cong  C^{\infty} (\mbox{$\bigwedge^3 T^*M$} ) \,
= \Omega^3, \notag
\end{equation}
and $C^{\infty}(PM)$ becomes a Riemannian manifold with respect to
the $L^2$ metric on $\Omega^3$ determined by $g_{\phi}$. Note the
diffeomorphism group Diff$(M)$ acts on $C^{\infty}(PM)$ as a group of
isometry.

\begin{definition}
A positive differential 3-form $\phi$ is a \textnormal{$G_2$ form} if
$\, d \phi = 0$ and $ \, d *_{\phi} \phi = 0$. A compact oriented manifold
$M$ is a \textnormal{$G_2$ manifold} if it admits a $G_2$ form.
\end{definition}

A $G_2$ metric $g_{\phi}$ is necessarily Ricci-flat. The
holonomy of the $G_2$ metric on a $G_2$ manifold is
isomorphic to a subgroup of $G_2$, and it is full $G_2$ whenever
the fundamental group $\pi_1(M)$ is finite, or equivalently the first de Rham
cohomology $H^1(M, R) = 0$.

We denote the space of $G_2$ forms by
$\mathfrak{\widehat{M}} \, \subset \, C^{\infty}(PM)$. Let
$\mathcal{D} \subset$ Diff$(M)$ be the group of diffeomorphisms of $M$
isotopic to the identity.
\begin{definition}
Let $M$ be a $G_2$ manifold. The moduli space of $G_2$ structures
(forms) is the quotient space
$\mathfrak{M} \, = \, \mathfrak{\widehat{M}} /  \mathcal{D}$.
\end{definition}
Given a $G_2$ form $\phi$, we denote its equivalence class by
$\langle \phi \rangle$.

It is known that $\mathfrak{M}$ is a smooth manifold of dimension
$b^3 =$ dim $H^3(M, R)$ [Br][Jo]. From the remark above, there exists
a
unique Riemannian metric on $\mathfrak{M}$ for which
$\mathfrak{\widehat{M}} \, \to \, \mathfrak{M}$ is a Riemannian
submersion. When $b^1 =$ dim $H^1(M, R) =0$, a connected component of
$\mathfrak{M}$ coincides with the Einstein deformation space of the
underlying $G_2$ metric, which has a real analytic structure [Ko].

From the definition of $\mathcal{D}$, the projection
$\pi^3 : \mathfrak{M} \to H^3(M, R)$ is well defined.
The image of a class is denoted by $\pi^3(\langle \phi \rangle) =
[\phi]$
for simplicity.

\noindent \emph{Remark}. $\;$ Let $f$ be a diffeomorphism of $M$. Then
\begin{equation}
*_{f^* \phi} f^* \phi \, = \, f^* ( *_{\phi} \phi ),
\end{equation}
and the cohomology map $\pi^4 : \mathfrak{M} \to H^4(M, R)$ by
$\pi^4(\langle \phi \rangle) = [*_{\phi} \phi]$ is also well defined.
In fact, we may take the dual definition of the moduli space of $G_2$
structures as the equivalence class of \emph{definite} 4-forms $\psi$
that satisfies $d \psi =0$ and $d  *_{\psi} \psi =0$.
(3) shows these two definitions "commute" with the projection maps
$\pi^3$ and $\pi^4$. Note however a definite 4-form does not
determine an orientation.

Let $\phi \in \mathfrak{\widehat{M}}$. Then $\bigwedge^* T^*M$ admits
a $G_2$ invariant decomposition
\begin{align}
\mbox{$\bigwedge^2$} T^*M \, &= \, \mbox{$\bigwedge^2_7 \,  \oplus \,
\bigwedge^2_{14}$} \notag \\
\mbox{$\bigwedge^3$} T^*M \, &= \, \mbox{$\bigwedge^3_1 \, \oplus
\,   \bigwedge^3_7 \, \oplus \, \bigwedge^3_{27}$} \notag
\end{align}
where, for instance,
\begin{equation}
*_{\phi} \phi \, \wedge \, (\, \mbox{$\bigwedge^3_7 \, \oplus \,
\bigwedge^3_{27}$} \,) \, = \, 0. \notag
\end{equation}

In particular, for any $X \in \Omega^3$, there exists a unique
quadratic form $h_X$ and a vector field $v_X$ such that
\begin{equation}
X \, = \, h_X \, \cdot \phi + v_X \, \lrcorner \, *_{\phi} \phi \notag
\end{equation}
where $h_X \, \cdot \phi$ means the action of $h_X$ as a derivation.

We denote $\Omega^p_k \, = \, C^{\infty}(\bigwedge^p_k)$ and write
$X = X_1 + X_7 + X_{27}$ for a given $X \in \Omega^3$ for its
irreducible components. Since Hodge Laplacian commutes with this
decomposition, [Ch], let $H^p_k$ be the corresponding decomposition
of the cohomology $H^p(M, R)$.

$\Omega^3 \cong T_{\phi} C^{\infty}(PM)$ admits a $L^2$ decomposition
along $\mathfrak{\widehat{M}}$ as follows [Jo].
\begin{align}
T_{\phi} \mathfrak{\widehat{M}} &= \, \{ \, X \in \Omega^3 \, | \; \;
dX=0, \quad d *_{\phi} ( \frac{4}{3}X_1 + X_7 - X_{27}) = 0 \, \} \\
V_{\phi} &= \{ \, d \, \Omega^2_7 \, \} \notag \\
H_{\phi} &= V_{\phi}^{\perp} \, \cap T_{\phi}
\mathfrak{\widehat{M}}
            \notag \\
         &= \{  \, X \in \Omega^3 \, | \; \; d X=0,  \quad  d
*_{\phi} X=0 \, \} \notag \\
N_{\phi} &= (T_{\phi} \mathfrak{\widehat{M}})^{\perp}.
\notag
\end{align}
Here $*_{\phi}(\frac{4}{3}X_1 + X_7 - X_{27})$ is the derivative of
the map $\phi \to *_{\phi} \phi$.
$H_{\phi}$ and $V_{\phi}$ represent the horizontal and vertical
subspaces of $T_{\phi} \mathfrak{\widehat{M}}$ with respect to the
submersion $\mathfrak{\widehat{M}} \, \to \, \mathfrak{M}$. The
orthogonal projection map from $T_{\phi} C^{\infty}(PM)$ to these
subspaces will be denoted by $\Pi^V_{\phi}$, $\Pi^H_{\phi}$,
and $\Pi^N_{\phi}$ respectively.

\section{Horizontal geodesics on $\mathfrak{\widehat{M}} \to
\mathfrak{M}$}
Let $\{ w^1, w^2, \, .. \, w^7 \}$ be a local coframe on $M$, and
\begin{equation}
X = X_{ijk} \, w^i \wedge w^j \wedge w^k \in
T_{\phi} C^{\infty}(PM)  \cong \Omega^3, \notag
\end{equation}
where $X_{ijk}$ is skew symmetric in all of its indices.
Fix $\phi \in  \mathfrak{\widehat{M}}$ so that $\langle w^i, w^j
\rangle_{g_{\phi}} = g_{\phi}^{ij}$. The $L^2$ inner product on
$T_{\phi} C^{\infty}(PM)$ is defined by
\begin{equation}
\ll X, Y \gg_{\phi} = \, 6 \,   \int_M \;   X_{ijk} \,
Y_{i'j'k'}
  \, g_{\phi}^{ii'} \, g_{\phi}^{jj'} \, g_{\phi}^{kk'} \,
dvol_{g_{\phi}}.
\end{equation}
Note $\ll \phi, \phi \gg_{\phi} = 7 \, Vol_{\phi}(M)$.

Let $\nabla$ be the Levi-Civita connection on $C^{\infty}(PM)$. Since
tangent vectors on $C^{\infty}(PM)$ can be identified with $\Omega^3$
valued functions,
\begin{equation}
\nabla_X \, Y \, = \, X(Y) + D_X\, Y
\end{equation}
where $X(Y)$ is the directional derivative of $Y$ as an $\Omega^3$
valued function, and $D_X\, Y$
is the covariant derivative of $Y$ considered as a translation
invariant vector field.

$D_X\, Y$ can be computed by the following lemma [Br].
\begin{lemma}
Let $\, Z \, = \, h_Z \, \cdot \,  \phi + v_Z \, \lrcorner \,
*_{\phi}
\phi \, \in \Omega^3$ and consider a curve $\, \phi_t = \phi + t \, Z
+
O(t^2) \in \widehat{\mathfrak{M}}$.
Then $g_{\phi_t} = g_{\phi} + t \, 2 h_Z + O(t^2)$.
\end{lemma}

\begin{proposition}
Let $\, X, \, Y, \, Z \in \Omega^3$ be translation invariant vector
fields on  $C^{\infty}(PM)$. Then
\begin{align}
 2 \ll D_X Y, \, Z \gg_{\phi} &= -2 \ll h_X \cdot Y + h_Y \cdot X ,
\, Z \gg_{\phi}
  + 2 \ll Y, h_Z \cdot X \gg_{\phi} \notag \\
  + & \int_M \; \mbox{tr}_{g_{\phi}}  (h_X)
 \langle Y, \, Z \rangle_{g_{\phi}} \; \mbox{dvol}_{g_{\phi}}
    + \int_M \; \mbox{tr}_{g_{\phi}}  (h_Y)                  \langle
X, \, Z  \rangle_{g_{\phi}} \; \mbox{dvol}_{g_{\phi}} \notag \\
    - & \int_M \; \mbox{tr}_{g_{\phi}}  (h_Z)
\langle X, \, Y  \rangle_{g_{\phi}} \; \mbox{dvol}_{g_{\phi}}. \notag
\end{align}
\end{proposition}

\noindent \emph{Proof.} \; Differentiate (5) and use the fact $[X, Y]
= [Y, Z] = [Z, X] =0$. $\square$

\vspace{1pc}

Let $\{ \gamma_t \} \subset \mathfrak{M}$ be a geodesic and $\phi_t
\in \mathfrak{\widehat{M}}$ be one of its horizontal lifts, which is
also
a geodesic in $\mathfrak{\widehat{M}}$.  As a curve in
$C^{\infty}(PM)$,
\begin{align}
\Pi^N_t \, (\phi^{'}_t) &= 0 \\
\Pi^N_t \, (\nabla_{\phi^{'}_t}\, {\phi^{'}_t} )
   &=    \nabla_{\phi^{'}_t}\, {\phi^{'}_t},
\end{align}
where $\Pi^N_t = \Pi^N_{\phi_t}$, and
\begin{equation}
\nabla_{\phi^{'}_t} \, {\phi^{'}_t} = \phi^{''}_t + D_{\phi^{'}_t} \,
\phi^{'}_t.
\end{equation}

From (8),
\begin{equation}
\phi^{''}_0 = \Pi^N_0 \, (\phi^{''}_0) -
            \Pi^{V+H}_0 \, ( D_{\phi^{'}_0} \, \phi^{'}_0).
\end{equation}
Differentiating (7),
 \begin{equation}
\frac{d}{dt} \, \Pi^N_t \, (\phi^{'}_0) \, |_{t=0} \; + \Pi^N_0 \,
(\phi^{''}_0) = 0,
\end{equation}
and we obtain
\begin{equation}
 \phi^{''}_0 =  - \frac{d}{dt} \, \Pi^N_t \, (\phi^{'}_0) \, |_{t=0}
\;
            - \; \Pi^{V+H}_0 \,  (D_{\phi^{'}_0} \, \phi^{'}_0).
\end{equation}

We record the following for later application.

\begin{lemma}
Let $\phi_t \in \mathfrak{\widehat{M}}$ be a curve and put
$\psi_t = *_{\phi_t} \phi_t$. Let $X \in V_{\phi_0} \oplus
H_{\phi_0}$.
Then
\begin{equation}
\int_M \; \psi_0 \wedge \frac{d}{dt} \,
\Pi^N_t ( X ) |_{t=0} \; = 0,
\notag
\end{equation}
\end{lemma}
where $X$ is considered as a translation invariant vector field along
$\phi_{t}$.

\emph{Proof.} \; Since $X$ is closed, $\Pi^N_t ( X )$ is
closed all $t$, and form the decomposition (4), $\Pi^N_t ( X )$ is in
fact exact.  $\square$

\section{Morse functions on   $ \mathfrak{M}_1 $ }

Let $\mathfrak{\widehat{M}}_1 \subset \mathfrak{\widehat{M}}$ be the
set of $G_2$ forms of volume 1, and $\, \mathfrak{M}_1 \subset
\mathfrak{M}$ its image under the projection $\mathfrak{\widehat{M}}
\to \mathfrak{M}$. Then $\, \mathfrak{M} \cong
\mathfrak{M}_1 \, \times R^+$, and $\mathfrak{M}_1$ is an embedded
hypersurface of $\mathfrak{M}$. In this section, we show each nonzero
element in $ H^4(M, R)$ becomes a Morse function when composed with
the cohomology map  $\mathfrak{M}_1 \, \to H^3(M, R)$.

Let $\phi \in \mathfrak{\widehat{M}}_1$.
Then the gradient of the volume function on $\mathfrak{\widehat{M}}$
at $\phi$ is $\frac{1}{3} \phi \in H_{\phi}$, (19). Let
$\nu_{\phi} = \frac{1}{\sqrt{7}} \phi \in H_{\phi}$ be the unit
normal to $\mathfrak{\widehat{M}}_1 \subset \mathfrak{\widehat{M}}$.
The second fundamental form of $\mathfrak{\widehat{M}}_1$ is
\begin{equation}
\mbox{II} \, = \, - \ll \textnormal{d} \nu, \textnormal{d} \phi \gg =
- \frac{1}{\sqrt{7}} \, \ll \textnormal{d} \phi, \textnormal{d} \phi
\gg, \notag
\end{equation}
and $\mathfrak{\widehat{M}}_1 \subset \mathfrak{\widehat{M}}$ is a
umbilic hypersurface. Note that d$\phi$ in the expression above is
not the exterior derivative,  but it is a tautological 1-form that
represents the infinitesimal displacement of $\phi$ in
$\mathfrak{\widehat{M}}_1$.

We wish to introduce a class of functions on $\mathfrak{M}$. Let
$\beta \in H^4(M, R)$ be a nonzero class, and define
\begin{equation}
F^{\beta} ( \langle \phi  \rangle) \,  =  \, \beta([\phi])
\end{equation}
where $\beta$ is regarded as an element of $H^3(M, R)^*$. Let
$F^{\beta}_1$ denote its restriction to $\mathfrak{M}_1$, and
$\mbox{Crit}(F^{\beta}_1)$ its critical point set.

\noindent \emph{Remark}. $\;$ From the \emph{Remark} below
\textbf{Definition
2}, we may also define
\begin{equation}
G^{\alpha}(\langle \phi \rangle) = \alpha([*_{\phi} \phi] ), \notag
\end{equation}
for $\alpha \in H^3(M, R)$.

\begin{proposition}
\begin{equation}
\textnormal{Crit}(F^{\beta}_1) \; = \; \{ \,  \langle \phi  \rangle
\, |  \, \,
\beta = c^{\beta}_{ \langle \phi  \rangle} \, [*_{\phi} \phi] \, \}
\notag
\end{equation}
for some constant $c^{\beta}_{ \langle \phi \rangle}$.
\end{proposition}
Note $7 \, c^{\beta}_{ \langle \phi  \rangle} = F^{\beta}_1 ( \langle
\phi  \rangle) \neq 0 $.

\vspace{1pc}
\noindent \emph{Proof.} \; From the description of the unit normal
$\nu$ above, $\phi$ is a critical point whenever $\beta$ annihilates
$H^3_7 \oplus H^3_{27}$. The proposition follows for $H^3(M, R)^* =
H^4(M, R)$.  $\square$

\vspace{1pc}

Let $ \langle \phi_0  \rangle \in \mbox{Crit}(F^{\beta}_1)$ and put
$\psi_0 = *_{\phi_0} \phi_0$. Then one finds
\begin{align}
\nabla^2 \, F^{\beta}_1 \, |_{ \langle \phi_0  \rangle} \, &= \,
\nabla^2 \, F^{\beta}   \, |_{ \langle \phi_0  \rangle} \, +
\frac{\partial}{\partial \nu} \,  F^{\beta}|_{ \langle \phi_0
\rangle} \, \mbox{II}_{ \langle \phi_0  \rangle}  \notag \\
&= \,   \nabla^2 \, F^{\beta} \, |_{ \langle \phi_0  \rangle} \,
 - c^{\beta}_{ \langle \phi_0  \rangle} \ll \mbox{d}\phi,
\mbox{d}\phi \gg_{ \langle \phi_0  \rangle}.
\end{align}

Now, let $X \in H_{\phi_0}$ be a horizontal lift of a tangent vector
$x \in T_{\langle \phi _0 \rangle} \mathfrak{M}_1$ and let $\phi_t
\subset \mathfrak{\widehat{M}}$ be a horizontal geodesic with
$\phi^{'}_0 = X$. From (12) and \textbf{Lemma 2},
\begin{align}
 \nabla^2 \, F^{\beta}(x,x) &= c^{\beta}_{\langle \phi _0 \rangle}
\int_M \, \psi_0 \wedge \phi^{''}_0  \notag \\
   &= - c^{\beta}_{\langle \phi _0 \rangle} \, \ll \phi_0, \, D_X X
\gg_{\phi_0},
\end{align}
for $\phi_0 \in H_{\phi_0}$. \textbf{Proposition 1} then gives,
\begin{align}
\ll \phi_0, \, D_X X \gg_{\phi_0}
 &= -2 \ll h_X \cdot X, \, \phi_0 \gg_{\phi_0}
 + \ll X, X \gg_{\phi_0} \\
&\; - \frac{7}{6} \ll X, X \gg_{\phi_0}
 \;  \mbox{since} \;  h_{\phi_0} = \frac{1}{3} \, g_{\phi_0} \;
\mbox{and} \; X_1=0 \notag \\
 &= -\frac{1}{6} \ll X, X \gg_{\phi_0}
   -2 \ll h_X \cdot X, \phi_0 \gg_{\phi_0}  \notag \\
 &= -\frac{1}{6} \ll X, X \gg_{\phi_0}
      -2 \ll X_{27},X_{27} \gg_{\phi_0}. \notag
\end{align}

\begin{theorem} Let $F^{\beta}_1$ be a function on $\mathfrak{M}_1$
defined in (13) with $\beta \in H^4(M, R)$. Then at a critical point
$\langle \phi_0 \rangle \in \mathfrak{M}_1$,
\begin{equation}
\nabla^2 \, F^{\beta}_1 \, |_{\langle \phi _0 \rangle}(x,x) \, = \,
-c^{\beta}_{\langle \phi _0 \rangle} ( \, \frac{5}{6} \ll X_7, X_7
\gg_{\phi_0} - \frac{7}{6} \ll X_{27}, X_{27} \gg_{\phi_0} \, ) \notag
\end{equation}
where $X \in H_{\phi_0}$ is the horizontal lift of $x \in T_{\langle
\phi _0 \rangle} \mathfrak{M}_1$, and $\beta = c^{\beta}_{\langle
\phi _0 \rangle} \, [*_{\phi_0} \phi_0]$ with  $7 \,
c^{\beta}_{\langle \phi _0 \rangle} = F^{\beta}_1 (\langle \phi _0
\rangle)$. $F^{\beta}_1$ is a Morse function on $\mathfrak{M}_1$ for
any nonzero $\beta \in H^4(M, R)$.
\end{theorem}
Note that in case  $H^1(M, R) = 0$, every critical point is either a
positive local minimum or a negative local maximum.

\vspace{1pc}

It is known that  $H^4_1$ component of the Pontryagin class
$p_1(M)$ is  $p_1(M)_1$ $=$ $- \frac{1}{56 \pi^2} \Arrowvert \,
 R_{g_{\phi}} \Arrowvert^2  \, [*_{\phi} \phi]$ for any $G_2$ form
$\phi \in \mathfrak{M}_{1}$ , where $\Arrowvert \, R_{g_{\phi}} \Arrowvert^2$
is the $L^2$ norm of the curvature tensor of the associated metric $g_{\phi}$.

\begin{corollary}
Let $M$ be a compact $G_2$ manifold with $p_1(M) \neq 0$. Let
$\mathfrak{M}_1$ be the moduli space of $G_2$ structures of volume 1.
Then each connected component of $\mathfrak{M}_1$ is noncompact.
\end{corollary}
\emph{Proof.} \;  $p_1(M) \neq 0$ implies $F^{p_1(M)}_1 <
0$ on $\mathfrak{M}_1$. The corollary follows from maximum
principle.  $\square$.

\begin{corollary}
Let $M$ be a compact $G_2$ manifold with  $b^1 =0$ and $b^3 = 2$.
Suppose $\phi_1$ and $\phi_2$ are isotopic $G_2$ forms, i.e., they
can be connected through $G_2$ forms, and that  $[\phi_2] = \lambda
[\phi_1]$ or $[*_{\phi_2}\phi_2] = \lambda^{\frac{4}{3}}
[*_{\phi_1}\phi_1]$ for some constant $\lambda > 0$. Then there exits a
diffeomorphism $f$ of $M$ isotopic to the identity such that
$\phi_2 = \lambda \, f^*(\phi_1)$. In particular,  each
$\pi^3 : \mathfrak{M} \to H^3(M, R)$ and  $\pi^4 : \mathfrak{M} \to
H^4(M, R)$ is one to one on every connected components of
$\, \mathfrak{M}$.
\end{corollary}
\emph{Proof.} \; Assume $\phi_1, \phi_2 \, \in
\mathfrak{M}_1$. Take $\beta =[*_{\phi_1}\phi_1]$. Since $b^1=0$ and
$b^3=2$, $F^{\beta}_1$ is a
positive function on the connected component of $\mathfrak{M}_1$ containing
$\langle \phi _1 \rangle$  with a unique critical point
$\langle \phi_1 \rangle$.

If  $[*_{\phi_2}\phi_2] = \lambda^{\frac{4}{3}} [*_{\phi_1}\phi_1]$,
it follows from considering $G^{\alpha}_1$ with $\alpha = [\phi_1]$.  $\square$

\vspace{1pc}

In case $b^1 = 0$, the image $[\mathfrak{M}_1] \subset H^3(M, R)$
displays the following geometric properties.

1. $[\mathfrak{M}_1] \subset H^3(M, R)$ is an immersed orientable
locally convex hypersurface that is transversal to the radial
direction.

2. Restriction to $[\mathfrak{M}_1]$ of every linear(height)
function in $H^4(M, R)$ is a Morse function with only positive local
minima or negative local maxima.

3. Take $\beta = p_1(M) \neq 0$, the first Pontryagin class. Then
$F^{\beta}_1(F^{\beta})$ is negative on
$\mathfrak{M}_1(\mathfrak{M})$ [Br].

4. Similarly, take $\beta= \sigma^2$ for any nonzero
$\sigma \in H^2(M, R)$. Then $F^{\beta}_1(F^{\beta})$ is negative on
$\mathfrak{M}_1(\mathfrak{M})$ [Br].

\vspace{1pc}
In view of \textbf{Theorem 1} and the properties listed above, we
propose the following conjecture on the compactness of $G_2$
structures of volume 1 with bounded $L^2$ norm of curvature.

\vspace{1pc}
\noindent \emph{Conjecture.} \, \,  Let $\{\phi_n\}$ be a sequence of
$G_2$ forms of volume 1 in a connected component of the moduli space
$\mathfrak{M}_1$  of $G_2$ forms on a compact manifold $M$ with $b^{1}=0$.
Suppose $\{ - p_1(M)([\phi_n)]) \} = \{ \frac{1}{56 \pi^2} \Arrowvert
\, R_{g_{\phi_n}} \Arrowvert^2 \}$ is bounded from above.
Then there exists a subsequence $\{ \phi_{n_k} \}$, a sequence of
diffeomorphisms $f_k$ isotopic to the identity,  and a $G_2$ form
$\phi$ of volume 1 such that $f_k^*  \phi_{n_k} \, \to \, \phi$
in $C^1$.
\vspace{1pc}

Suppose this conjecture is true, and consider the  function
$F^{\beta}_1$ with $\beta = -p_1(M)$. From the  $C^1$ convergence,
\begin{align}
[f_k^*  \phi_{n_k}] \, &\to \, [\phi] \notag \\
[f_k^* *_{\phi_{n_k}}  \phi_{n_k}] \, &\to \, [*_{\phi} \phi], \notag
\end{align}
and the conjecture implies the set
\begin{equation}
 (F^{\beta}_1)^c = \{ \langle \phi \rangle \in \mathfrak{M}_1 \, | \,
F^{\beta}_1(\langle \phi \rangle) \leq c \} \notag
\end{equation}
is compact for any constant $c$.
Since every critical points of $F^{\beta}_1$ is a positive local
minimum, a result form Morse theory shows that there exists a unique
critical point $\langle \phi \rangle$ in each connected component of
$\mathfrak{M}_1$, and that every connected component of
$\mathfrak{M}_1$ is contractible. By \textbf{Theorem 1},
$[*_{\phi} \phi]$ must be a constant multiple of $- p_1(M)$.
Since $\mathfrak{M} \cong \mathfrak{M}_1 \times R^+$,  this implies
every connected component of $\mathfrak{M}$ is contractible and
contains a unique element $\langle \phi \rangle$ such that $[*_{\phi}
\phi] = -p_1(M)$. Such $\langle \phi \rangle$ is in fact unique up to
diffeomorphism, for Pontryagin class is a diffeomorphism invariant.

\section{Torelli map}
\begin{definition}
Let $M$ be a $G_2$ manifold,  and let $\mathfrak{M}$ be the moduli
space of $G_2$ structures on $M$. Torelli map
$\, \pi : \mathfrak{M} \, \to  H^3(M, R) \oplus H^4(M, R)$ is defined
by
\begin{equation}
\pi(\langle \phi \rangle) = (\pi^3(\langle \phi \rangle),
\pi^4(\langle \phi \rangle)) = ( [\phi], [*_{\phi} \phi]). \notag
\end{equation}
\end{definition}

For the period map of Riemann surfaces of fixed genus,
the image is known to be an open set of an irreducible analytic
subset [Gr].
An essential characterization of the image of the period map,
inversion problem, is also known via automorphic forms [AM]. The
purpose of this section is to describe the image of the Torelli map
$\pi(\mathfrak{M})$, or the projections $\pi^3(\mathfrak{M}_1)$ and
$\pi^4(\mathfrak{M}_1)$ in more detail. A general idea is that any
invariant of $\pi(\mathfrak{M})$ as a Lagrangian submanifold, or the
invariants of  $\pi^3(\mathfrak{M}_1)$ as an affine hypersurface, is
an invariant of $\mathfrak{M}$ or $\mathfrak{M}_1$.

Set $l=b^3$, and let  $\{ \, \alpha_1, \alpha_2, \, .. \, \alpha_l \,
\}$ be a basis of   $H^3(M, R)$  and $\{ \, \beta^1, \beta^2, \, ..
\, \beta^l \, \}$ be the dual basis of $H^4(M, R)$  so that
$ \beta^A(\alpha_B)=\delta^A_B$. Let
$x^t=(x^1, x^2, \, ..\, x^l)$ and $y^t=(y_1, y_2, \, .. \, y_l)$ be
the  coordinate functions on $H^3(M, R)$ and
$H^4(M, R)$ with respect to these basis. Then $\sum_A dx^{A} \wedge
dy_{A}$ is the canonical symplectic form, and
$G_0 = \sum_A dx^{A} \, dy_{A}$ is the canonical metric of
signature $(l, l)$ on  $H^3(M, R) \oplus H^4(M, R)$.

From the definition,
\begin{align}
\pi^3(\langle \phi \rangle) &= \mbox{$\sum_A$} \,  x^{A}(\langle \phi
\rangle) \alpha_A \\
\pi^4(\langle \phi \rangle) &= \mbox{$\sum_A$} \,  y_{A}(\langle \phi
\rangle) \beta^A, \notag
\end{align}
where
\begin{align}
x^{A}(\langle \phi \rangle)  &= \beta^A([\phi]) \\
y_{A}(\langle \phi \rangle)  &=\alpha_A([*_{\phi} \phi]). \notag
\end{align}
Each $x$ and $y$ is then a local diffeomorphism from $\mathfrak{M}$
to $R^l$ [Jo].

Define a function $U :  \mathfrak{M} \to R$ by
\begin{equation}
U(\langle \phi \rangle) = \mbox{$\sum_A$} x^A \, y_{A} = 7 \,
Vol_{\phi}(M) = \int_M \, \phi \wedge *_{\phi} \phi. \notag
\end{equation}
Then,
\begin{align}
dU &= \sum_A x^A \, dy_{A} + y_{A} dx^A  \\
    &= \int_M \phi \wedge \,  d\pi^4 + d\pi^3 \wedge \, *_{\phi} \phi.
\notag \\
    &= \frac{7}{3} \int_M d\pi^3 \wedge \,  *_{\phi} \phi. \notag \\
    &= \frac{7}{3}  \sum_A  y_{A} dx^A,  \notag
\end{align}
where the third equality follows from (4).
Hence
\begin{equation}
3 \sum_A x^A dy_{A} \, = \, 4 \sum_A y_{A} dx^A,
\end{equation}
and the Torelli map  $\pi$ is a Lagrangian immersion [Jo].

Let $p= (p_{AB}) = (p_{BA})$ be the unique $Gl(l, R)$ valued
function on $\mathfrak{M}$ such that
\begin{equation}
dy_{A} = \mbox{$\sum_B$} \, p_{AB} dx^B. \notag
\end{equation}
From (20), we have
\begin{equation}
4 \, y \, = \, 3 \, p \, x.
\end{equation}
Thus the invertible symmetric matrix function $p$ transforms one
period function to the other.

Upon a change of basis $x^* = A^{-1} x$, $y^*= A^t y$ for $A \in
Gl(l, R)$, $p$ becomes $p^* = A^t \, p \, A$, and $dx^t \, p \, dx =
dx^t \, dy$ is a well defined quadratic form on $\mathfrak{M}$, which
is  the metric induced from $G_0$ via  Torelli map $\pi$.
\begin{proposition}
The signature of $\pi^*G_0$ is
$(1+b^1, b^3_{27} = b^3-b^1-1)$.
If $\, b^1=0$, $\pi^*G_0$ is the negative of the $L^2$ metric when
restricted to $\mathfrak{M}_1 \subset \mathfrak{M}$.
\end{proposition}
\emph{Proof.} \; Let $X \in T_{\phi}
\mathfrak{\widehat{M}}$. Then, with an abuse of notation,
\begin{align}
\pi^*G_0(X,X) &= dx^t(X) \, dy(X) \notag \\
              &= \frac{4}{3} \ll X_1, X_1 \gg_{\phi} +
\ll X_7, X_7 \gg_{\phi}
-\ll X_{27}, X_{27} \gg_{\phi} \notag
\end{align}
by (4). $\square$

\vspace{1pc}

Let $\{ \xi_{A}(\phi) \}_{1}^{l}$ be a basis of $H^{3}(M, R)$ where
$[\xi_{A}(\phi)]=\alpha_{A}$ and $\xi_{A}(\phi)$ is harmonic with
respect to $g_{\phi}$. Define $m_{AB}(\langle \phi \rangle)$ by
\begin{equation}
[*_{\phi} \xi_{A}(\phi) ] = \sum_{B} \, m_{AB} \, \beta^{B}, \notag
\end{equation}
or equivalently
\begin{equation}
m_{AB}(\langle \phi \rangle) =
\int_{M} \, *_{\phi} \xi_{A}(\phi) \, \wedge
\xi_{B}(\phi),
\end{equation}
and
\begin{equation}
    m_{AB} \, dx^{A} \, dx^{B} \notag
\end{equation}
is the $L^{2}$ metric on $\mathfrak{M}$.
Suppose $b^{1}=0$. From \textbf{Proposition 3}, we get
\begin{align}
p_{AB} &= -m_{AB}+ \frac{7}{3U} \, y_{A} \, y_{B} \\
p^{AB} &= -m^{AB}+ \frac{7}{4U} \, x^{A} \, x^{B}, \notag
\end{align}
where $(p^{AB})$ and $(m^{AB})$ are the inverse matrices of
$(p_{AB})$ and $(m_{AB})$ respectively.

\vspace{1pc}

We now turn our attention to the hypersurface
 $\Sigma = [\mathfrak{M}_V] \subset H^3(M, R)$, where
$\mathfrak{M}_V \subset \mathfrak{M}$ is the set of $G_2$ structures
of volume $V$.
 For definiteness, let us assume $\{ \alpha_1, \alpha_2, \, .. \,
\alpha_l \}$ is a basis of $H^3(M, Z)$ modulo torsion, choose an
orientation for $H^3(M, R)$, and consider the properties of $\Sigma$
that are invariant under the linear change of basis by $Sl(l, R)$
[Ch2][NS]. We agree on the index range
$1 \, \leq i, \, j \, \leq l-1$ and  $1 \, \leq A, \, B \, \leq l$.

Let
\begin{align} %24
\vec{x} &= \sum_A x^A \, \alpha_A  \\
\vec{y} &= \sum_A y_A \, \beta^A  \notag
\end{align}
be the immersions defined by (17), (18). Then since
$\mathfrak{M}_V$ is defined by the equation $U=7 \, V$,
\begin{equation}
%25
\sum_A y_A dx^A=\sum_A x^A dy_A=0
\end{equation}
on $\mathfrak{M}_V$, and we may write
\begin{align} %26
d\vec{x} &= \mbox{$\sum_i$} dx^i \, \alpha_i + dx^l \, \alpha_l \\
   &= \mbox{$\sum_i$} dx^i  (\alpha_i - \frac{y_i}{y_l} \alpha_l) +
(dx^l + \frac{1}{y_l} \mbox{$\sum_i$}y_i dx^i) \, \alpha_l \notag \\
   &= \mbox{$\sum_i$} \omega^i \, e_i + \omega^l \, e_l
\notag
\end{align}
where
\begin{align} %27
e_i &= \alpha_i - \frac{y_i}{y_l} \alpha_l \\
e_l &= \alpha_l \notag
\end{align}
and
\begin{align} %28
\omega^i &= dx^i \\
\omega^l &= dx^l + \frac{1}{y_l} \mbox{$\sum_i$}y_i dx^i = 0. \notag
\end{align}
Following the general theory of moving frames, put
\begin{equation} %29
de_A = \mbox{$\sum_B$}\omega_A^B \, e_B. \notag
\end{equation}
Differentiating (27), we get
\begin{align} %30
de_i &\equiv \omega^l_i \, e_l \quad \mod \, \, \, e_1, e_2, \, .. \,
e_{l-1} \\
     &\equiv -d (\frac{y_i}{y_l}) e_l, \notag
\end{align}
and
\begin{align}
\omega^l_i &=-d (\frac{y_i}{y_l}) \\
      &=  \, \mbox{$\sum_j$} \, h_{ij} \omega^j \notag
\end{align}
where
\begin{align}
h_{ij} &= h_{ji} \notag \\
       &= -\frac{p_{ij}}{y_l} -\frac{p_{ll}}{y_l^3} y_i y_j +
\frac{y_{li} y_j+y_{lj}y_i}{y^2_l}. \notag
\end{align}
By (28), the second fundamental form $II$ of $\Sigma$
 \begin{align} %32
II &= \mbox{$\sum_{ij}$} h_{ij} \, \omega^i \, \omega^j \notag \\
   &= -\frac{1}{y_l} \mbox{$\sum_{AB}$}p_{AB} dx^A \, dx^B,
\end{align}
and it is definite whenever $b^{1}=0$ by \textbf{Proposition 3}.

Set $H=$ det$(h_{ij}) \neq 0$. Then the normalized second fundamental
form
\begin{equation} %33
\hat{II} = |H|^{-\frac{1}{l+1}} \, II
\end{equation}
is an affine invariant called \emph{Blaschke metric}.

\noindent \emph{Remark}. If $b^1=0$, one may assume
$H > 0$ on $\Sigma$.

Since
\begin{equation}
\omega^l_1 \wedge \omega^l_2 \, .. \, \wedge \omega^l_{l-1} = H \,
\omega^1 \wedge \omega^2 \, .. \, \wedge \omega^{l-1}, \notag
\end{equation}
a computation using (30) gives,
\begin{align}%34 35
\omega^l_1 \wedge \omega^l_2 \, .. \, \wedge \omega^l_{l-1} &=
\frac{1}{y^l_l} \, (\sum_A y_A \, \frac{\partial}{\partial y_A} \,
\lrcorner \, dy_1 \wedge dy_2 \, .. \, \wedge dy_l) \\
&=\frac{1}{y^l_l} \, \frac{3}{4} \, \mbox{det}(p) (\sum_A x^A \,
\frac{\partial}{\partial x^A} \, \lrcorner \, dx^1 \wedge dx^2 \, ..
\, \wedge dx^l)
\end{align}
for $\sum_A y_A \, \frac{\partial}{\partial y_A} = \frac{3}{4} \sum_A x^A \, \frac{\partial}{\partial x^A}$ by (21). By
expanding (34) and using $U=7V$,
\begin{equation}%36
H=\frac{21}{4} V \mbox{det}(p) (-1)^{l-1} y_l^{-(l+1)},
\end{equation}
and
\begin{equation}%37
\hat{II}= -(\frac{21}{4}V)^{-\frac{1}{l+1}}
|\mbox{det}(p)|^{-\frac{1}{l+1}} \, \mbox{$\sum_{AB}$}p_{AB} dx^A
dx^B.
\end{equation}
Here we assume $\{ \, \alpha_1, \alpha_2, \, .. \, \alpha_l \, \}$ is
the right orientation for
$H^3(M, R)$.

Another affine invariant called \emph{affine normal} $\xi^x$ is
defined by
\begin{equation} %38
\xi^x = |H|^{\frac{1}{l+1}} \, (e_l + \mbox{$\sum_{i}$} t^i \, e_i),
\end{equation}
where $t^i$'s are uniquely determined by the equation
\begin{equation} %39
\frac{1}{l+1} \, d \, log |H| + \mbox{$\sum_{i}$} t^i \omega^l_i \, =
0 = \lambda \mbox{$\sum$}_A x^A dy_A,
\end{equation}
for some nonzero $\lambda$.

Set
\begin{equation}%40
q^A = \frac{1}{(l+1)\mbox{det}(p)}\frac{\partial
\mbox{det}(p)}{\partial y_A}.
\notag
\end{equation}
Then from (35) and (30), $\lambda= \frac{1}{7V}(-1+\mbox{$\sum$}_A
y_A q^A)$, and we get
\begin{equation}%41
t^i = y_l (q^i - x^i \lambda).
\end{equation}
The affine normal $\xi^x$ is now determined to be
\begin{align}%42
\xi^x &= |H|^{\frac{1}{l+1}} \, (e_l + \mbox{$\sum_{i}$} t^i \, e_i),
\\
 &=  |H|^{\frac{1}{l+1}} \, \frac{y_l}{7V} \Big(
(\mbox{$\sum$}_A x^A \alpha_A) + (7Vq^A - x^A(\mbox{$\sum$}_B y_B
q^B))\alpha_A \Big) \notag \\
 &= \frac{1}{7V} (\frac{21}{4}V)^{\frac{1}{l+1}}
|\mbox{det}(p)|^{\frac{1}{l+1}} \Big(
(\mbox{$\sum$}_A x^A \alpha_A) + (7Vq^A - x^A(\mbox{$\sum$}_B y_B
q^B))\alpha_A \Big). \notag
\end{align}

A hypersurface in an affine space is  an \emph{affine sphere} if all
the normal lines pass through a fixed point(finite or infinite)
called the center. In case it is locally convex, it is called
elliptic or hyperbolic depending on whether the center is on the
convex or concave side.

\begin{theorem}
The locus of projection $\pi^3(\mathfrak{M}_1) \subset \, H^3(M, R)$
is an affine sphere centered at the origin if det$(p)$ is constant
along $\mathfrak{M}_1$. In case $b^{1}=0$, this is equivalent to
$det(m)$ being constant on $\mathfrak{M}_1$, and
$\pi^3(\mathfrak{M}_1)$ is a (locally convex) hyperbolic affine sphere.
\end{theorem}
\noindent \emph{Remark}. $\;$ $det(m)$ is constant on $\mathfrak{M}_1$
if and only if the volume of the torus
$ \, H^{3}(M, R)/H^3(M, Z) $ is constant on $\mathfrak{M}_1$.

\noindent \emph{Proof.} $\;$ $\xi^x$ is proportional to $\vec{x}$ if
$q^A = \mu x^A$ for all $A$ for some function $\mu$ on $\mathfrak{M}_1$.
By (40) and (25), this implies $d$ det$(p) =0$ along $\mathfrak{M}_{1}$.

In case $b^{1}=0$,
\begin{equation}
    \sum_{AB} p^{AB} dp_{AB} \equiv \sum_{AB} m^{AB} dm_{AB}
    \; \, \mod \; \, dU \notag
\end{equation}
from (23). $\pi^3(\mathfrak{M}_1)$ is locally
convex by \textbf{Proposition 3} and hyperbolic by \textbf{Theorem 1}.  $\square$

\vspace{1pc}

Similar computation for the conormal hypersurface
$\Sigma^{'} = \pi^4(\mathfrak{M}_1) \subset H^4(M, R)$
defined by $\vec{y}$ in (24) gives the affine normal
\begin{equation}
\xi_y = \frac{1}{7V}
(\frac{3}{28V})^{\frac{1}{l+1}}|\mbox{det}(p^{-1})|^{\frac{1}{l+1}}
\Big( (\mbox{$\sum$}_A y_A \beta^A) + (-7Vq_A + y_A(\mbox{$\sum$}_B
x^B q_B))\beta^A \Big),
\end{equation}
with
\begin{align}
q_A &= \frac{1}{(l+1)\mbox{det}(p)}\frac{\partial
\mbox{det}(p)}{\partial x^A} \notag \\
    &= \sum_B p_{AB} q^B. \notag
\end{align}
Note the functions
\begin{align}
\xi_y(\vec{x}) &=(\frac{3}{28V})^{\frac{1}{l+1}}
|\mbox{det}(p^{-1})|^{\frac{1}{l+1}} \notag \\
\xi^x(\vec{y}) &= (\frac{21V}{4})^{\frac{1}{l+1}}
|\mbox{det}(p)|^{\frac{1}{l+1}} \notag \\
\xi_y(\xi^x) &= \frac{1}{7V} (\frac{9}{16})^{\frac{1}{l+1}} \Big( 1 +
(\mbox{$\sum$}_A x^A q_A)(\mbox{$\sum$}_B y_B q^B) - 7 V \,
\mbox{$\sum$}_A q^A q_A \Big). \notag
\end{align}
are independent of the choice of the orientation for $H^3(M, R)$.
From (19),
\begin{align}
\sum_A x^A \frac{\partial \mbox{det}(p)}{\partial x^A} &= \frac{l}{3}
\mbox{det}(p)   \notag    \\
\sum_A y_A \frac{\partial \mbox{det}(p)}{\partial y_A} &= \frac{l}{4}
\mbox{det}(p). \notag
\end{align}
Hence
\begin{equation}
\xi_y(\xi^x) = \frac{1}{7V}(\frac{9}{16})^{\frac{1}{l+1}}\Big( 1 +
\frac{l^2}{12(l+1)^2} -
7 V\, \mbox{$\sum$}_A q^A q_A \Big). \notag
\end{equation}

\newpage

\noindent \begin{center}
\textbf{\large{References}}
\end{center}

\noindent
[AM] A. Andreotti, A. Mayer, On period relations on abelian integral
for algebraic curves, Ann. Scuola  Norm. Sup. Pisa CI. Sci. 21
(1967), 189 - 238

\noindent
[Be] A. Besse, \emph{Einstein manifolds}, Springer-Verlag (1987)

\noindent
[Br] R. Bryant, $G_2$ manifolds, seminar lectures, Duke University
(1999)

\noindent
[Ch1] S. S. Chern,  On a Generalization of K\"ahler Geometry, S. S.
Chern Selected Papers I, 227 - 245

\noindent
[Ch2] S. S. Chern, Affine minimal hypersurfaces, S. S. Chern Selected
Papers III, 425 - 438

\noindent
[Gr] P. Giffiths, Deformations of complex structure,
Proc. Sym. Pure Math. 15, 251 - 274 (1970)

\noindent
[Hi] N. Hitchin,   The moduli space of special Lagrangian
submanifolds,
Ann. Scuola  Norm. Sup. Pisa CI. Sci. (4), 25(3-4) (1997), 503 - 515

\noindent
[Jo] D. Joyce,  \emph{Compact manifolds with special holonomy},
Oxford University Press (2000)

\noindent
[Ko] N. Koiso,  Einstein metrics and complex structures, Invent.
Math. \textbf{73} (1983), 71 - 106

\noindent
[NS] K. Nomizu, T. Sasaki,  \emph{Affine differential geometry},
Cambridge University Press (1994)

\noindent
[Tr] A. Tromba, \emph{Teichm\"uller theory in Riemannian geometry},
Lecture Notes in Math. ETH, (1992)

\vspace{2pc}

\noindent Sung Ho Wang\\
Department of Mathematics\\
Postech\\
Pohang, Korea 790-784\\
\noindent \emph{email}: \textbf{wang@postech.ac.kr}

\end{document}